\documentclass[twocolumn]{autart}

\usepackage{amsmath,amssymb}
\usepackage{graphicx}
\usepackage[longnamesfirst,sort]{natbib}% Citation support using natbib.sty
\usepackage{subcaption}
\usepackage{xcolor}

\newcommand{\RR}{\mathbb{R}}

\newtheorem{Assumption}{Assumption}
\newtheorem{Remark}{Remark}
\newtheorem{Lemma}{Lemma}
\newtheorem{Theorem}{Theorem}
\newtheorem{Problem}{Problem}

\newenvironment{Proof}{\noindent{\em Proof:\/}}{\hfill $\Box$\par}

\begin{document}

\begin{frontmatter}

\title{Distributed exponential state estimation of linear systems over jointly connected switching networks\thanksref{footnoteinfo}}

\thanks[footnoteinfo]{
%This work has been supported
%in part by National Natural Science
%Foundation of China under the grant No.\ 61973260,
%and in part by the Research Grants Council of the Hong Kong Special Administrative Region under the grants Nos.\ 14201420 and PDFS2122-4S04.
Corresponding author: J.\ Huang.
Tel.\ +852-39438473. Fax +852-39436002.
}

\author[]{Tao Liu}\ead{tliu2@mae.cuhk.edu.hk},
\author[]{Jie Huang}\ead{jhuang@mae.cuhk.edu.hk}

\address{Department of Mechanical and Automation Engineering, The Chinese University of Hong Kong, Shatin, N.T., Hong Kong}


\begin{abstract}
Recently, the distributed state estimation problem for continuous-time linear systems over jointly connected switching networks was solved.
It was shown that the estimation errors will asymptotically converge to the origin by using the generalized Barbalat's Lemma.
This paper further studies the same problem with two new features.
First, the asymptotic convergence is strengthened to the exponential convergence.
This strengthened result not only offers a guaranteed convergence rate,
but also renders the error system total stability and thus is able to withstand small disturbances.
Second, the coupling gains of our local observers can be distinct and thus offers greater design flexibility,
while the coupling gains in the existing result were required to be identical.
These two new features are achieved by establishing exponential stability for two classes of linear time-varying systems,
which may have other applications.

\end{abstract}

\end{frontmatter}

%%%%%%%%%%%%%%%%%%%%%%%%%%%%%%%%%%%%%%%%%%%%%%%%%%%%%%%%%%%%%%%%%%%%%%%%%%%%%%%%%%%%%%%%%%%%%%%%%%%%%%%%%%%%%%%%%%%%%%%%%%%%%%%%%%%%%%%%%%%%%%%%%%%%%%%%%%%%%%%%%%%%%%%%%%%%%%%%%%%%%%%%%%%%%%%%%%
\section{Introduction}
Designing observers to estimate the state of a given plant is one of the fundamental problems in modern control theory \citep{Luenberger71}.
Recently, driven by the rapid development in both theory and applications of sensor networks and multi-agent systems,
the distributed state estimation problem has attracted considerable attentions.
The distributed state estimation problem aims to design a network of local observers,
in which each local observer can only access partial information of the given plant.
One typical scenario from which this problem arises is
a large-scale complex system, such as a smart grid or an industrial process,
where the system is monitored by a group of spatially distributed sensors
that transmitting information over a communication network.

The distributed state estimation problem was first studied over static networks
by \cite{Acikmese14}, \cite{ParkMartins17}, and \cite{MitraSundaram18} for discrete-time linear systems
and by \cite{KimShim16} and \cite{WangMorse18} for continuous-time linear systems.
In particular, \cite{Acikmese14} first proposed a two-time-scale design of discrete-time local observers to solve the problem.
Each local observer in \cite{Acikmese14} is assisted by a consensus filter
that operates multiple times between every successive pair of time instants of the system.
A first single-time-scale design of discrete-time local observers was given by \cite{ParkMartins17}.
They presented a parameterized class of discrete-time local observers with certain augmented states,
and they established a necessary and sufficient condition for the existence of local observers
by casting the problem into a stabilization problem for some augmented system via fully decentralized output feedback control.
This idea was further explored by \cite{WangMorse18} to deal with the continuous-time case.
In addition, \cite{WangMorse18} developed an approach via the decentralized control theory in \cite{CorfmatAnderson76}
to freely assign the spectrum of the overall observer so that arbitrarily fast convergence rate can be achieved.
\cite{MitraSundaram18} studied a similar problem to that of \cite{ParkMartins17}
by first performing a so-called multi-sensor observable canonical decomposition on the plant.
Then, based on the block lower triangular form resulting from this decomposition,
local observers without state augmentation were devised to solve the problem.
It should be noted that the first Luenberger-type local observers
were proposed by \cite{KimShim16} for distributed state estimation of continuous-time linear systems,
which inspired several interesting extensions in, e.g., \cite{HanTrentelman19} and \cite{WangLiu20}.
Also, \cite{KimShim16} further refined their design of Luenberger-type local observers to the current form in \cite{KimShim20}.
%Some other variants of the distributed state estimation problem were studied in some recent papers, say, \cite{Slim19}, \cite{Xu20}, and \cite{WangSu23}.

The existing literature on the distributed state estimation problem mainly focuses on static and connected networks.
The first extension to switching networks was made by \cite{WangLiu20}.
Nevertheless, the switching networks in \cite{WangLiu20} were required to be strongly connected for all time,
while, in practice, disconnectedness of a network may be caused by intermittent communication link failures or environmental changes.
Thus, it is more interesting and challenging to further study the distributed estimation problem over jointly connected switching networks,
which can be disconnected at every time instant. Indeed, a significant advance was made recently by \cite{ZhangLu21ACC}.
By using a common Lyapunov function approach in conjunction with the generalized Barbalat's Lemma established in \cite{SuHuang12Full},
\cite{ZhangLu21ACC} showed that the time-varying version of the local observers of \cite{KimShim20} is able to
asymptotically estimate the state of a class of neutrally stable linear systems over jointly connected switching networks.

The objective of this paper is to strengthen the result of \cite{ZhangLu21ACC}
from the asymptotic convergence to the exponential convergence under the same conditions as those in \cite{ZhangLu21ACC}.
The strengthened result offers at least two advantages.
First, the exponential convergence result leads to the guaranteed convergence rate,
which is much desired in practice.
Second, since, for linear systems, exponential stability is equivalent to uniform asymptotic stability \citep{Rugh96},
the strengthened result implies that the error system is totally stable and hence is able to withstand small disturbances \citep{Slotine91}.
For this purpose, instead of using the generalized Barbalat's Lemma as adopted in \cite{ZhangLu21ACC},
we have to come up with a completely new approach that makes use of the classical uniformly complete observability concept.
We need to first establish the uniformly complete observability for a class of linear time-varying systems.
This result then leads to two exponential stability results for two classes of linear time-varying systems
in Lemmas \ref{Lemma-eta-stable} and \ref{Lemma-zeta-stable}, respectively.
As a result of Lemmas \ref{Lemma-eta-stable} and \ref{Lemma-zeta-stable},
we conclude exponential stability for the overall estimation error system.
Moreover, as a byproduct, we show that the coupling gains of our local observers can be distinct,
while the coupling gains in \cite{ZhangLu21ACC} were required to be identical.
This new feature offers greater flexibility in the design of local observers.

\medskip

\noindent
\emph{Notation.}
$\RR$ denotes the set of real numbers.
For $X_{i}\in \RR^{m_{i}\times n}, i=1,\ldots,N$,
$\text{col}\left(X_{1}, \ldots, X_{N}\right)=\left[
                                               \begin{array}{ccc}
                                                 X_{1}^{T} & \cdots & X_{N}^{T} \\
                                               \end{array}
                                             \right]^{T}
$.
For a matrix $A\in \RR^{m\times n}$, we denote its kernel by
$\text{ker}(A)=\{x\in \RR^{n} : Ax=0\}$
and its range by $\text{im}(A)=\{y \in \RR^{m}: y=Ax\ \text{for some}\ x\in \RR^{n}\}$.
$\mathbf{1}_{N}$ denotes the $N$ dimensional column vector whose entries are all $1$.
$\mathbf{0}$ denotes a zero matrix with conformable dimensions.
$\otimes$ denotes the Kronecker product of matrices.
For any $x\in \RR^{n}$, $\|x\|$ denotes the Euclidean norm of $x$.
$\mathbf{e}$ denotes the base of the natural logarithm.
$\lambda_{\max}(Q)$ and $\lambda_{\min}(Q)$ denote the largest and the smallest eigenvalue of
a real symmetric matrix $Q$, respectively.

\noindent
\emph{Terminology.}
We call a time function $\sigma:[0,\infty)\mapsto \mathcal{P}:=\{1,2,\ldots,n_{0}\}$,
a piecewise constant switching signal,
if there exists a sequence $\{t_{j}: j=0,1,2,\ldots\}$ satisfying
$t_{0}=0$ and $t_{j+1}-t_{j}\ge \tau, j=0,1,2,\ldots$, for some $\tau>0$,
such that for all $t\in [t_{j},t_{j+1}), \sigma(t)=p$ for some $p\in \mathcal{P}$.
Then $\mathcal{P}$ is called the switching index set,
$\{t_{j}: j=0,1,2,\ldots\}$ are called switching instants,
and $\tau$ is called the dwell time.

Given a linear time-varying system
\begin{equation}
    \dot{z}(t)=A(t)z(t), \quad y(t)=C(t)z(t), \quad t\ge 0 \tag{$\star$}
\end{equation}
where $z(t)$ is the state, $y(t)$ is the output,
and $A(t), C(t)$ are time-varying matrices of conformable dimensions.
Let $\Phi_{z}(t,t^{*}), t\ge t^{*} \ge 0$ denote the state transition matrix of system ($\star$),
and for any $t^{*}+T \ge t^{*} \ge 0$, let the observability Gramian of system ($\star$) be denoted by
\begin{equation*}
    G_{z}(t^{*},t^{*}+T):=\int_{t^{*}}^{t^{*}+T}\Phi_{z}(t,t^{*})^{T}C(t)^{T}C(t)\Phi_{z}(t,t^{*})\, d\, t.
\end{equation*}
Then, system ($\star$) is said to be uniformly completely observable (UCO, see pp.\ 35 of \cite{Sastry89})
if there exist $T_{\text{o}}>0$ and $0<\alpha_{1}\le \alpha_{2}$ such that
\begin{equation*}
   \alpha_{1}I \le G_{z}(t^{*},t^{*}+T_{\text{o}}) \le \alpha_{2}I,\quad \forall\ t^{*} \ge 0.
\end{equation*}

\section{Problem Formulation}\label{Section_Problem_Formulation}
Consider the following linear time-invariant system:
\begin{equation}\label{eq-x-Ax}
    \dot{x}(t)=Ax(t), \qquad y(t)=Cx(t), \qquad t\ge 0
\end{equation}
where $x(t)\in \RR^{n}$ is the state,
$y(t)\in \RR^{m}$ is the output,
$A\in \RR^{n\times n}$ is the system matrix,
and $C\in \RR^{m \times n}$ is the output matrix.

Suppose there is a network of $N$ agents and each agent can only measure partial information of system \eqref{eq-x-Ax} as follows:
\begin{equation}\label{eq-y-i}
    y_{i}(t)=C_{i}x(t), \qquad i=1,\ldots,N
\end{equation}
where, for $i=1,\ldots,N$,
$y_{i}(t)\in \RR^{m_{i}}$ is the partial output of system \eqref{eq-x-Ax} measured by the $i$th agent,
$\sum_{i=1}^{N}m_{i}=m$, and
$C_{i}\in \RR^{m_{i}\times n}$ are such that $\text{col}(C_{1}, \ldots, C_{N})=C$.

The topology of the communication network for these $N$ agents is described by
an undirected switching graph\footnote{See Appendix A for a summary of notation on graph.}
$\mathcal{G}_{\sigma(t)}=(\mathcal{V},\mathcal{E}_{\sigma(t)})$,
where $\sigma(t)$ is the piecewise constant switching signal,
$\mathcal{V}=\{1,\ldots,N\}$, and $\mathcal{E}_{\sigma(t)}\subseteq \mathcal{V}\times \mathcal{V}$.
Specifically, each node in $\mathcal{V}$ corresponds to an agent in the network,
and $(i,j)\in \mathcal{E}_{\sigma(t)}$ if and only if the $i$th agent can exchange information with the $j$th agent at the time instant $t$.

Now we are ready to describe the distributed \emph{exponential} state estimation problem.

\begin{Problem}\label{Problem}
  Given system \eqref{eq-x-Ax}, the local measurements \eqref{eq-y-i},
  and the switching graph $\mathcal{G}_{\sigma(t)}$,
  design for each agent, a local observer of the following form:
  \begin{equation}\label{eq-hat-x-i-general}
    \dot{\hat{x}}_{i}(t)=\mathbf{f}_{i}\left(\hat{x}_{i}(t),y_{i}(t),\{\hat{x}_{j}(t)-\hat{x}_{i}(t): j \in \mathcal{N}_{i}(t)\}\right),
     t\ge 0
  \end{equation}
  where $\hat{x}_{i}(t) \in \RR^{n}$ is the observer state, $\mathbf{f}_{i}(\cdot)$ is some linear function,
  and $\mathcal{N}_{i}(t)=\{j \in \mathcal{V} : (j,i) \in \mathcal{E}_{\sigma(t)}\}$,
  such that, for any initial conditions $x(0)\in \RR^{n}$ and $\hat{x}_{i}(0)\in \RR^{n}, i=1,\ldots,N$,
  the solutions of systems \eqref{eq-x-Ax} and \eqref{eq-hat-x-i-general} satisfy
  $\lim_{t\to\infty}(\hat{x}_{i}(t)-x(t))=0, i=1,\ldots,N$, exponentially.
\end{Problem}

\begin{Remark}
  It is interesting to compare the above problem with the distributed observer design problem studied in, say,
  \cite{SuHuang12Cyber}, \cite{CaiHuang16}, and \cite{LTHuang19}.
  The problem here divides the output $y(t)$ into $N$ components and then the aim is to design  $N$
  local observers over a communication network of $N$ nodes to estimate the state of system \eqref{eq-x-Ax}.
  On the other hand, the distributed observer design problem involves $(N+1)$
  subsystems consisting of $N$ followers and one leader which is system \eqref{eq-x-Ax}.
  The communication among these $(N+1)$ subsystems is governed by a graph of $(N+1)$ nodes.
 Assuming, at each time instant, only a subset of the followers can access the output $y(t)$ of system \eqref{eq-x-Ax},
 one needs to design $N$ local observers  over a communication network of $(N+1)$ nodes to estimate the state of system \eqref{eq-x-Ax}.
\end{Remark}

In order to solve Problem \ref{Problem}, we make the following three assumptions
on the system matrix $A$, the pair $(C,A)$, and the switching graph $\mathcal{G}_{\sigma(t)}$, respectively.

\begin{Assumption}\label{Ass-A-neutral}
The matrix $A$ is neutrally stable, i.e., all the eigenvalues of $A$ are semi-simple with zero real parts.
\end{Assumption}

Under Assumption \ref{Ass-A-neutral}, there exist a nonsingular matrix $P \in \RR^{n \times n}$ such that
  the matrix $A$ is similar to a skew-symmetric matrix, i.e.,
  \begin{equation*}
    \bar{A}:=P^{-1}AP \qquad \text{and} \qquad \bar{A}^{T}=-\bar{A}.
  \end{equation*}
For simplicity of presentation, in what follows, we assume that the matrix $A$ is skew-symmetric.
The design procedure for the general case where the matrix $A$ is not skew-symmetric
will be outlined in Remark \ref{Remark-design-procedure}.

\begin{Assumption}\label{Ass-observable}
  The pair $(C,A)$ is observable.
\end{Assumption}

%\begin{Remark}\label{Remark-observable}
%Under Assumption \ref{Ass-observable},
%the observability matrix $\mathcal{O}$ of the pair $(C,A)$ is of full column rank, i.e.,
%\begin{equation*}\label{eq-rank-O}
%    \text{rank}\left(\mathcal{O}:=\left[
%                                    \begin{array}{c}
%                                      C \\
%                                      CA \\
%                                      \vdots \\
%                                      CA^{n-1} \\
%                                    \end{array}
%                                  \right]
%     \right)=n.
%\end{equation*}
%In other words, $\text{ker}(\mathcal{O})=\{0\}$.
%\end{Remark}

\begin{Assumption}\label{Ass-jointly-connected}
  There exists a subsequence $\{t_{j_{k}}: k=0,1,2,\ldots\}$
  of the switching instants $\{t_{j}: j=0,1,2,\ldots\}$
  with $t_{j_{0}}=0$ and $ \tau \le t_{j_{k+1}}-t_{j_{k}} \le T_{\text{c}}$
  for some $T_{\text{c}}>0$, such that the union graph
  $\bigcup_{t_{j}\in [t_{j_{k}}, t_{j_{k+1}})}\mathcal{G}_{\sigma(t_{j})}$ is connected.
\end{Assumption}

\begin{Remark}\label{Remark-jointly-connected}
Let $\mathcal{L}_{\sigma(t)}$ be the Laplacian of the undirected switching graph $\mathcal{G}_{\sigma(t)}$.
Then, $\mathcal{L}_{\sigma(t)}$ is symmetric and positive semi-definite for all $t \ge 0$.
Moreover, $\sum_{r=j_{k}}^{j_{k+1}-1}\mathcal{L}_{\sigma(t_{r})}$ is a Laplacian associated with the union graph
$\bigcup_{t_{j}\in [t_{j_{k}}, t_{j_{k+1}})}\mathcal{G}_{\sigma(t_{j})}$.
Then, under Assumption \ref{Ass-jointly-connected},
the matrix $\sum_{r=j_{k}}^{j_{k+1}-1}\mathcal{L}_{\sigma(t_{r})}$ has exactly one zero eigenvalue
and its null space is spanned by the single vector $\mathbf{1}_{N}$ \citep{Lin06,SuHuang12Full}.
\end{Remark}

\begin{Remark}
Assumptions \ref{Ass-A-neutral} to \ref{Ass-jointly-connected} were also used in \cite{ZhangLu21ACC}.
More specifically, Assumption \ref{Ass-jointly-connected} is called the jointly connected condition \citep{Jadbabaie03, SuHuang12Full}
or the uniformly connected condition \citep{Lin06}.
It is perhaps the mildest condition for a graph as it allows the graph to be disconnected at any time instant.
Hence, to compensate for the weak connectivity of the graph under Assumption \ref{Ass-jointly-connected},
one has to impose some stability constraint on the system to be observed, namely, Assumption \ref{Ass-A-neutral}.
It is worth pointing out that Assumption \ref{Ass-A-neutral} has been a standard assumption
in dealing with consensus problems for general linear multi-agent systems over jointly connected switching networks \citep{SuHuang12Full}.
It should also be noted that the two classes of most popular multi-agent systems
including single-integrator dynamics with $A=\mathbf{0}$ and oscillator dynamics
with $A =\text{block diag}\left\{\left[
                                   \begin{array}{cc}
                                     0 & 1 \\
                                     -\alpha & 0 \\
                                   \end{array}
                                 \right],\cdots, \left[
                                   \begin{array}{cc}
                                     0 & 1 \\
                                     -\alpha & 0 \\
                                   \end{array}
                                 \right] \right\}$
and $\alpha>0$ satisfy Assumption \ref{Ass-A-neutral}
\citep{Olfati-Saber04, Ren08}.
\end{Remark}

\section{Main results}\label{Section_Main_Results}
To begin with, for the $i$th agent, $i=1,\ldots,N$,
given the pair $(C_{i}, A)$,
there exists a nonnegative integer $\nu_{i}$ such that
\begin{equation*}
    \text{rank}\left(\mathcal{O}_{i}:=\left[
                                    \begin{array}{c}
                                      C_{i} \\
                                      C_{i}A \\
                                      \vdots \\
                                      C_{i}A^{n-1} \\
                                    \end{array}
                                  \right]
     \right)=n-\nu_{i}.
\end{equation*}
Hence, $\text{ker}(\mathcal{O}_{i})\subseteq \RR^{n}$ is the $\nu_{i}$-dimensional unobservable subspace of $(C_{i},A)$.

%\begin{Remark}\label{Remark-joint-observability}
%  It can be readily seen that $\bigcap_{i=1}^{N}\text{ker}(\mathcal{O}_{i})=\text{ker}(\mathcal{O})$.
%  Thus, under Assumption \ref{Ass-observable}, by Remark \ref{Remark-observable},
%  if there exists some $z \in \RR^{n}$ such that $z \in \text{ker}(\mathcal{O}_{i})$ for $i=1,\ldots,N$,
%  then $z=0$.
%\end{Remark}

Let $U_{i} \in \RR^{n \times \nu_{i}}$ be a matrix
whose columns form an orthonormal basis of $\text{ker}(\mathcal{O}_{i})$,
namely,
\begin{equation}\label{eq-U-i-def}
    \text{im}(U_{i})=\text{ker}(\mathcal{O}_{i})
\end{equation}
and let $D_{i} \in \RR^{n \times (n-\nu_{i})}$ be a matrix whose columns form an orthonormal basis of $\text{im}(\mathcal{O}_{i}^{T})$,
namely,
\begin{equation*}
    \text{im}(D_{i})=\text{im}(\mathcal{O}_{i}^{T}).
\end{equation*}
Then, we can define an orthogonal matrix $T_{i}$ as follows:
\begin{equation}\label{eq-T-i-def}
    T_{i}:=\left[
             \begin{array}{cc}
               D_{i} & U_{i} \\
             \end{array}
           \right]\in \RR^{n \times n}, \qquad T_{i}^{T}T_{i}=I_{n}.
\end{equation}

\begin{Lemma}\label{Lemma-transformation}
Under Assumption  \ref{Ass-A-neutral},
for $i=1,\ldots,N$, we have
 \begin{equation}\label{eq-T-i-A-T-i}
    T_{i}^{T}AT_{i}=\left[
                          \begin{array}{cc}
                            A_{i\text{o}} & \mathbf{0} \\
                            \mathbf{0} & A_{i\bar{\text{o}}} \\
                          \end{array}
                        \right], \qquad C_{i}T_{i}=\left[
                                                    \begin{array}{cc}
                                                      C_{i\text{o}} & \mathbf{0} \\
                                                    \end{array}
                                                  \right]
\end{equation}
where $A_{i\text{o}}\in \RR^{(n-\nu_{i})\times (n-\nu_{i})}$,
$A_{i\bar{\text{o}}}\in \RR^{\nu_{i}\times \nu_{i}}$,
and $C_{i\text{o}} \in \RR^{m_{i}\times (n-\nu_{i})}$.
Moreover,
\begin{enumerate}
  \item[(i)] $(C_{i\text{o}},  A_{i\text{o}})$ is observable;

  \item[(ii)] $A_{i\bar{\text{o}}}$ is skew-symmetric;

  \item[(iii)] $A U_{i}=U_{i} A_{i\bar{\text{o}}}$.
\end{enumerate}
\end{Lemma}

\begin{Remark}
  Equation \eqref{eq-T-i-A-T-i} is a direct result of the standard Kalman decomposition and from which, Part (i) follows.
  Part (ii) is obvious since,  under Assumption \ref{Ass-A-neutral}, $A$ is assumed to be skew-symmetric.
  Part (iii) is a consequence of the fact that the subspace $\text{ker}(\mathcal{O}_{i})$ is $A$-invariant.
  It is noted that the decomposition \eqref{eq-T-i-A-T-i} was also performed in Lemma 3.2 of \cite{ZhangLu21ACC}
  together with a detailed construction of the matrix $T_{i}$.
\end{Remark}

Next, to introduce our first stability result,
let $\nu=\sum_{i=1}^{N}\nu_{i}$ and
\begin{equation}\label{eq-U-block-def}
     U=\text{block diag}\{U_{1},U_{2},\ldots, U_{N}\} \in \RR^{Nn \times \nu}.
\end{equation}

\begin{Lemma}\label{Lemma-eta-stable}
  Consider the following linear switched system:
  \begin{equation}\label{eq-eta-system}
    \dot{\eta}(t)=- Q U^{T}\left(\mathcal{L}_{\sigma(t)} \otimes I_{n}\right)U\eta(t), \qquad t\ge 0
  \end{equation}
  where $\eta(t) \in \RR^{\nu}$ is the state
  and $Q \in \RR^{\nu\times \nu}$ is any symmetric and positive definite matrix.
  Under Assumptions \ref{Ass-observable} and \ref{Ass-jointly-connected},
  system \eqref{eq-eta-system} is exponentially stable.
\end{Lemma}
\emph{Proof:} See Appendix B.

\begin{Remark}
If $Q=\gamma I_{\nu}$ with $\gamma>0$,
then system \eqref{eq-eta-system} reduces to the following linear switched system:
 \begin{equation}\label{eq-eta-system-special}
    \dot{\eta}(t)=- \gamma U^{T}\left(\mathcal{L}_{\sigma(t)} \otimes I_{n}\right)U\eta(t), \qquad t\ge 0.
 \end{equation}
What makes Lemma \ref{Lemma-eta-stable} interesting is that,
under Assumptions \ref{Ass-observable} and \ref{Ass-jointly-connected},
system \eqref{eq-eta-system} or \eqref{eq-eta-system-special} is exponentially stable even though the following system:
 \begin{equation*}
    \dot{\eta}(t)=-\gamma \left(\mathcal{L}_{\sigma(t)} \otimes I_{n}\right) \eta(t), \qquad t\ge 0
  \end{equation*}
can only be marginally stable, no matter what condition is imposed on the graph $\mathcal{G}_{\sigma(t)}$.
\end{Remark}

Now we are ready to establish the stability result for another class of linear switched systems.
To state the next lemma, let
\begin{align}
    A_{\bar{\text{o}}}&=\text{block diag}\{A_{1\bar{\text{o}}}, A_{2\bar{\text{o}}},\ldots, A_{N\bar{\text{o}}}\} \in \RR^{\nu \times \nu} \label{eq-A-o-block-def}\\
\Gamma &= \text{block diag}\{\gamma_{1}I_{\nu_{1}},\gamma_{2}I_{\nu_{2}},\ldots, \gamma_{N}I_{\nu_{N}}\} \in \RR^{\nu \times \nu} \label{eq-Gamma-block-def}
\end{align}
where $A_{i\bar{\text{o}}} \in \RR^{\nu_{i} \times \nu_{i}}, i=1,\ldots,N$, are given in Lemma \ref{Lemma-transformation}
 and $\gamma_{i}>0, i=1,\ldots,N$.

\begin{Lemma}\label{Lemma-zeta-stable}
 Consider the following linear switched system:
  \begin{equation}\label{eq-zeta-system}
    \dot{\zeta}(t)=\left(A_{\bar{\text{o}}}-\Gamma U^{T}\left(\mathcal{L}_{\sigma(t)} \otimes I_{n}\right)U \right)\zeta(t),
    \quad t\ge 0
  \end{equation}
  where $\zeta(t) \in \RR^{\nu}$ is the state.
  Under Assumptions \ref{Ass-A-neutral} to \ref{Ass-jointly-connected},
  system \eqref{eq-zeta-system} is exponentially stable.
\end{Lemma}

\begin{Proof}
  Let $\eta(t)=\mathbf{e}^{-A_{\bar{\text{o}}}t}\zeta(t), t\ge 0$.
  Then, by noticing that $A$ is skew-symmetric, $A_{i\bar{\text{o}}}$ is also skew-symmetric as shown in Part (ii) of Lemma \ref{Lemma-transformation},
  and $AU_{i}=U_{i}A_{i\bar{\text{o}}}$ as shown in Part (iii) of Lemma \ref{Lemma-transformation},
  we can derive that
  \begin{align*}
    \dot{\eta}(t)&=-\mathbf{e}^{-A_{\bar{\text{o}}}t} A_{\bar{\text{o}}} \zeta(t)+\mathbf{e}^{-A_{\bar{\text{o}}}t} \dot{\zeta}(t) \notag \\
     &=-\mathbf{e}^{-A_{\bar{\text{o}}}t} \Gamma U^{T}\left(\mathcal{L}_{\sigma(t)} \otimes I_{n}\right)U \zeta(t) \notag \\
     &=- \Gamma\left[
                 \begin{array}{ccc}
                  \mathbf{e}^{A^{T}_{1\bar{\text{o}}}t}U_{1}^{T} &  &  \\
                    & \ddots &  \\
                    &  & \mathbf{e}^{A^{T}_{N\bar{\text{o}}}t}U_{N}^{T} \\
                 \end{array}
               \right]\left(\mathcal{L}_{\sigma(t)} \otimes I_{n}\right)U \zeta(t) \notag \\
     &=- \Gamma \left[
                 \begin{array}{ccc}
                   U_{1}^{T}\mathbf{e}^{A^{T}t} &  &  \\
                    & \ddots &  \\
                    &  & U_{N}^{T}\mathbf{e}^{A^{T}t}  \\
                 \end{array}
               \right]\left(\mathcal{L}_{\sigma(t)} \otimes I_{n}\right)U \zeta(t) \notag \\
     &=-\Gamma U^{T}\left(I_{N} \otimes \mathbf{e}^{-At}\right)\left(\mathcal{L}_{\sigma(t)} \otimes I_{n}\right)U \zeta(t) \notag\\
     &=-\Gamma U^{T}\left(\mathcal{L}_{\sigma(t)} \otimes I_{n}\right)\left(I_{N} \otimes \mathbf{e}^{-At}\right)U \zeta(t) \notag \\
     &=-\Gamma U^{T}\left(\mathcal{L}_{\sigma(t)} \otimes I_{n}\right)U\mathbf{e}^{-A_{\bar{\text{o}}}t}\zeta(t) \notag\\
     &=-\Gamma U^{T}\left(\mathcal{L}_{\sigma(t)} \otimes I_{n}\right)U\eta(t).
  \end{align*}
Then, the above $\eta$-system is in the form of system \eqref{eq-eta-system} with $Q=\Gamma$.
Thus, by Lemma \ref{Lemma-eta-stable}, for any $\eta(0)\in \RR^{\nu}$,
$\lim_{t\to\infty}\eta(t)=0$ exponentially.

Since, in addition,
$\eta(t)^{T}\eta(t)=\zeta(t)^{T}\mathbf{e}^{A_{\bar{\text{o}}}t}\mathbf{e}^{-A_{\bar{\text{o}}}t}\zeta(t)=\zeta(t)^{T}\zeta(t),
\forall\, t\ge 0$,
we can conclude that for any $\zeta(0)\in \RR^{\nu}$,
$\lim_{t\to\infty}\zeta(t)=0$ exponentially.
Therefore, system \eqref{eq-zeta-system} is exponentially stable,
and the proof is complete.
\end{Proof}

To present our main result,  for $i=1,\ldots,N$,  by extending the Luenberger-type local observer in \cite{KimShim20} to switching networks gives
the following local observer:
\begin{align}\label{eq-hat-x-i-specific}
    \dot{\hat{x}}_{i}(t)&=A\hat{x}_{i}(t)+L_{i}\left(y_{i}(t)-C_{i}\hat{x}_{i}(t)\right) \notag\\
    &\qquad \  + \gamma_{i} M_{i}\sum_{j=1}^{N}a_{ij}(t)\left(\hat{x}_{j}(t)-\hat{x}_{i}(t)\right), \quad t\ge 0
\end{align}
where $\hat{x}_{i}(t)\in \RR^{n}$ is the observer state,
$\gamma_{i}>0$ is the coupling gain of the local observer, which can be different from each other,
$a_{ij}(t)$ are entries of the weighted adjacency matrix $\mathcal{A}_{\sigma(t)}$ of the switching graph $\mathcal{G}_{\sigma(t)}$,
and $L_{i} \in \RR^{n \times m_{i}}$ and $M_{i}\in \RR^{n \times n}$ are the injection matrix and the weighting matrix, respectively.
Specifically,  $L_{i}$ and $M_{i}$ are designed as
\begin{equation}\label{eq-L-i-M-i-def}
   L_{i}=T_{i}\left[
                 \begin{array}{c}
                   L_{i\text{o}} \\
                   \mathbf{0} \\
                 \end{array}
               \right], \qquad M_{i}=T_{i} \left[
                                             \begin{array}{cc}
                                               \mathbf{0} & \mathbf{0} \\
                                               \mathbf{0} & I_{\nu_{i}} \\
                                             \end{array}
                                           \right]T_{i}^{T}
\end{equation}
where $T_{i}\in \RR^{n\times n}$ is the orthogonal matrix given by \eqref{eq-T-i-def},
and $L_{i\text{o}} \in \RR^{(n-\nu_{i})\times m_{i}}$ is such that
$(A_{i\text{o}}-L_{i\text{o}}C_{i\text{o}})$ is Hurwitz,
whose existence is guaranteed by Part (i) of Lemma \ref{Lemma-transformation}.

\begin{Theorem}\label{Theorem}
  Under Assumptions \ref{Ass-A-neutral} to \ref{Ass-jointly-connected},
Problem \ref{Problem} is solvable by designing for each agent $i, i=1,\ldots,N$, a local observer of the form \eqref{eq-hat-x-i-specific}.
\end{Theorem}

\begin{Proof}
For $i=1,\ldots,N$, let $e_{i}(t)=\hat{x}_{i}(t)-x(t)$ be the estimation error of the $i$th local observer.
Then, the error dynamics of $e_{i}(t)$ can be written as follows:
\begin{align}\label{eq-e-i-dynamics}
   \dot{e}_{i}(t)&=(A-L_{i}C_{i})e_{i}(t) + \gamma_{i} M_{i}\sum_{j=1}^{N}a_{ij}(t)\left(e_{j}(t)-e_{i}(t)\right) \notag\\
    &=(A-L_{i}C_{i})e_{i}(t) - \gamma_{i} M_{i}\sum_{j=1}^{N}l_{ij}(t)e_{j}(t)
\end{align}
where $l_{ij}(t)$ are entries of the Laplacian $\mathcal{L}_{\sigma(t)}$ of the switching graph $\mathcal{G}_{\sigma(t)}$.

Next, for $i=1,\ldots,N$, perform the following coordinate transformation on $e_{i}(t)$:
\begin{equation*}
    \left[
      \begin{array}{c}
       \xi_{i\text{o}}(t)  \\
       \xi_{i\bar{\text{o}}}(t) \\
      \end{array}
    \right]:= \left[
                \begin{array}{c}
                  D_{i}^{T} \\
                  U_{i}^{T} \\
                \end{array}
              \right]e_{i}(t)=T_{i}^{T}e_{i}(t)
\end{equation*}
where $\xi_{i\text{o}}(t) \in \RR^{n-\nu_{i}}$ and $\xi_{i\bar{\text{o}}}(t) \in \RR^{\nu_{i}}$.
Then, by using \eqref{eq-T-i-A-T-i} and \eqref{eq-L-i-M-i-def},
the error dynamics \eqref{eq-e-i-dynamics} can be further written as
\begin{align}\label{eq-xi-i-dynamics}
    \dot{\xi}_{i\text{o}}(t)&=\left(A_{i\text{o}}-L_{i\text{o}}C_{i\text{o}}\right)\xi_{i\text{o}}(t) \notag\\
    \dot{\xi}_{i\bar{\text{o}}}(t)&= A_{i\bar{\text{o}}}\xi_{i\bar{\text{o}}}(t)
    - \gamma_{i} U_{i}^{T}\sum_{j=1}^{N}l_{ij}(t)\left(D_{j} \xi_{j\text{o}}(t)+U_{j}\xi_{j\bar{\text{o}}}(t)\right).
\end{align}

Now, for the purpose of analyzing the stability of system \eqref{eq-xi-i-dynamics},
let
\begin{align*}
    \xi_{\text{o}}(t)&=\text{col}\left(\xi_{1\text{o}}(t), \xi_{2\text{o}}(t), \ldots, \xi_{N\text{o}}(t) \right) \\
    \xi_{\bar{\text{o}}}(t)&=\text{col}\left(\xi_{1\bar{\text{o}}}(t), \xi_{2\bar{\text{o}}}(t), \ldots, \xi_{N\bar{\text{o}}}(t) \right)
\end{align*}
and, in addition to \eqref{eq-U-block-def}, \eqref{eq-A-o-block-def}, and \eqref{eq-Gamma-block-def},
let $X=\text{block diag}\{X_{1}, X_{2}, \ldots, X_{N} \}$,
for $X_{i}=A_{i\text{o}}, L_{i\text{o}}, C_{i\text{o}}$, and $D_{i}, i=1,\ldots,N$.
Then,  the $N$ error dynamics in \eqref{eq-xi-i-dynamics} can be put into the following compact form:
\begin{align}\label{eq-xi-dynamics-compact}
    \dot{\xi}_{\text{o}}(t) &= \left(A_{\text{o}}-L_{\text{o}}C_{\text{o}}\right)\xi_{\text{o}}(t) \notag\\
    \dot{\xi}_{\bar{\text{o}}}(t) &= A_{\bar{\text{o}}}\xi_{\bar{\text{o}}}(t)
    - \Gamma U^{T}\left(\mathcal{L}_{\sigma(t)} \otimes I_{n}\right)\left(D \xi_{\text{o}}(t) +U\xi_{\bar{\text{o}}}(t) \right).
\end{align}

Since, by our design, the matrix $\left(A_{\text{o}}-L_{\text{o}}C_{\text{o}}\right)$ is Hurwitz,
$\xi_{\text{o}}$-subsystem of \eqref{eq-xi-dynamics-compact} is exponentially stable.
Then, we have $\lim_{t\to\infty}\Gamma U^{T}\left(\mathcal{L}_{\sigma(t)} \otimes I_{n}\right)D \xi_{\text{o}}(t)=0$ exponentially.
By Lemma 1 of \cite{LTHuang18},
for any initial condition $\xi_{\bar{\text{o}}}(0) \in \RR^{\nu}$,
the solution $\xi_{\bar{\text{o}}}(t)$ of $\xi_{\bar{\text{o}}}$-subsystem of \eqref{eq-xi-dynamics-compact} converges to zero exponentially
if the following linear switched system:
\begin{equation}\label{eq-xi-dynamics-compact-unperturbed}
\dot{\xi}_{\bar{\text{o}}}(t) = \left(A_{\bar{\text{o}}}- \Gamma U^{T}\left(\mathcal{L}_{\sigma(t)} \otimes I_{n}\right)U \right)\xi_{\bar{\text{o}}}(t)
\end{equation}
is exponentially stable, which, in fact, has been shown in Lemma \ref{Lemma-zeta-stable}.
Therefore, system \eqref{eq-xi-dynamics-compact} is exponentially stable and hence,
for any initial conditions $x(0)\in \RR^{n}$ and $\hat{x}_{i}(0)\in \RR^{n}, i=1,\ldots,N$,
$\lim_{t\to\infty}(\hat{x}_{i}(t)-x(t))=0, i=1,\ldots,N$, exponentially.
The proof is thus complete.
\end{Proof}

\begin{Remark}
 It is interesting to note that for static graphs considered in \cite{KimShim20},
 the switching Laplacian $\mathcal{L}_{\sigma(t)}$ reduces to a constant $\mathcal{L}$,
 and it was shown in Lemma 4 of \cite{KimShim20} that $U^{T}\left(\mathcal{L}\otimes I_{n} \right)U$ is positive definite.
 Thus, in order to stabilize the unobservable part of the error dynamics,
 which is in a time-invariant form of \eqref{eq-xi-dynamics-compact-unperturbed},
 it suffices to choose sufficiently large $\gamma_{i}>0, i=1,\ldots,N$,
 to make the matrix  $\left(A_{\bar{\text{o}}}-\Gamma U^{T}\left(\mathcal{L} \otimes I_{n}\right)U \right)$ Hurwitz.
 Nevertheless, for switching graphs considered in this paper,
 system \eqref{eq-xi-dynamics-compact-unperturbed} is time-varying,
and one has to prove exponential stability for \eqref{eq-xi-dynamics-compact-unperturbed} by a completely different approach.
 \end{Remark}

\begin{Remark}
For the special case of system \eqref{eq-xi-dynamics-compact} where $\Gamma=\gamma I_{\nu}$ and $\gamma>0$,
an asymptotic stability result was obtained in \cite{ZhangLu21ACC} by treating the two subsystems in \eqref{eq-xi-dynamics-compact} as a whole.
This approach calls for the construction of a common Lyapunov function for the switched system \eqref{eq-xi-dynamics-compact} and
the usage of the generalized Barbalat's Lemma in \cite{SuHuang12Full}.
In contrast, we have managed to establish exponential stability for system \eqref{eq-xi-dynamics-compact-unperturbed} by using the UCO concept,
which in turn concludes exponential stability for system  \eqref{eq-xi-dynamics-compact} due to its lower triangular structure.
It is also interesting to note that the common Lyapunov function proposed in \cite{ZhangLu21ACC} works
only if $\gamma_{1}=\gamma_{2}=\cdots=\gamma_{N}=\gamma>0$.
Thus, Theorem \ref{Theorem} also extends the main result of \cite{ZhangLu21ACC}
by allowing different local observers to have different coupling gains $\gamma_{i}, i=1,\ldots,N$.
\end{Remark}

%\begin{Remark}
%  Suppose the local measurements in \eqref{eq-y-i} are specified with
%  $C_{1}=C$ and $C_{i}=\mathbf{0}, i=2,\ldots,N$.
%  Then, the local observers in \eqref{eq-hat-x-i-specific} become
%  \begin{align}\label{eq-hat-x-i-special-case}
%    \dot{\hat{x}}_{1}(t)&=A\hat{x}_{1}(t)+L\left(y(t)-C\hat{x}_{1}(t)\right), \quad t\ge 0 \notag\\
%    \dot{\hat{x}}_{i}(t)&=A\hat{x}_{i}(t)+\gamma \sum_{j=1}^{N}a_{ij}(t)\left(\hat{x}_{j}(t)-\hat{x}_{i}(t)\right),\ i=2,\ldots,N
%\end{align}
%where $L\in \RR^{n \times m}$ is such that $(A-LC)$ is Hurwitz.
%By Theorem \ref{Theorem}, we also have $\lim_{t\to\infty}(\hat{x}_{i}(t)-x(t))=0, i=1,\ldots,N$, exponentially.
%It is noted that system \eqref{eq-hat-x-i-special-case} can be viewed as a simple variant of
%the distributed observer studied in \cite{SuHuang12Cyber}.
%In fact, it has been shown in \cite{SuHuang12Cyber} that, under much more relaxed assumptions,
%$\lim_{t\to\infty}(\hat{x}_{i}(t)-x(t))=0, i=1,\ldots,N$, exponentially still holds.
%Specifically, the switching graph $\mathcal{G}_{\sigma(t)}$ need not be undirected,
%Assumption \ref{Ass-A-neutral} can be relaxed to that $A$ has no eigenvalues with positive real parts,
%and Assumption \ref{Ass-jointly-connected} can be relaxed to that
%every node $i,i=2,\ldots,N$ is reachable from node $1$ in the union graph $\bigcup_{r=j_{k}}^{j_{k+1}-1}\mathcal{G}_{\sigma(t_{r})}$.
%\end{Remark}

\begin{Remark}\label{Remark-design-procedure}
  If the matrix $A$ is neutrally stable but not skew-symmetric, then the design of the local observer
  in \eqref{eq-hat-x-i-specific} can be carried out by the following procedure:
  \begin{enumerate}
    \item Find $P$ such that $P^{-1}AP$ is skew-symmetric.

    \item Find $U_{i}$ such that $\text{im}(U_{i})=\text{ker}(\mathcal{O}_{i}P)$.

    \item Find $D_{i}$ such that $\text{im}(D_{i})=\text{im}(P^{T}\mathcal{O}_{i}^{T})$.

  \item Form $T_{i}$ as in \eqref{eq-T-i-def} and perform the similarity transformation \eqref{eq-T-i-A-T-i} with $PT_{i}$ in place of $T_{i}$ to obtain $(C_{i\text{o}}, A_{i\text{o}})$.

    \item Design $L_{i\text{o}}$ such that $(A_{i\text{o}}-L_{i\text{o}}C_{i\text{o}})$ is Hurwitz and let
    $L_{i}=PT_{i}\left[
                 \begin{array}{c}
                   L_{i\text{o}} \\
                   \mathbf{0} \\
                 \end{array}
               \right]$, $M_{i}=PT_{i} \left[
                                             \begin{array}{cc}
                                               \mathbf{0} & \mathbf{0} \\
                                               \mathbf{0} & I_{\nu_{i}} \\
                                             \end{array}
                                           \right]T_{i}^{T}P^{-1}$.
  \end{enumerate}
\end{Remark}

\section{A numerical example}\label{Section_Example}
In this section, we use a modified example of Example 1 in \cite{KimShim20} to illustrate our design of the local observers over a jointly connected  switching network.

\begin{figure}
\centering
  \includegraphics[scale=0.19]{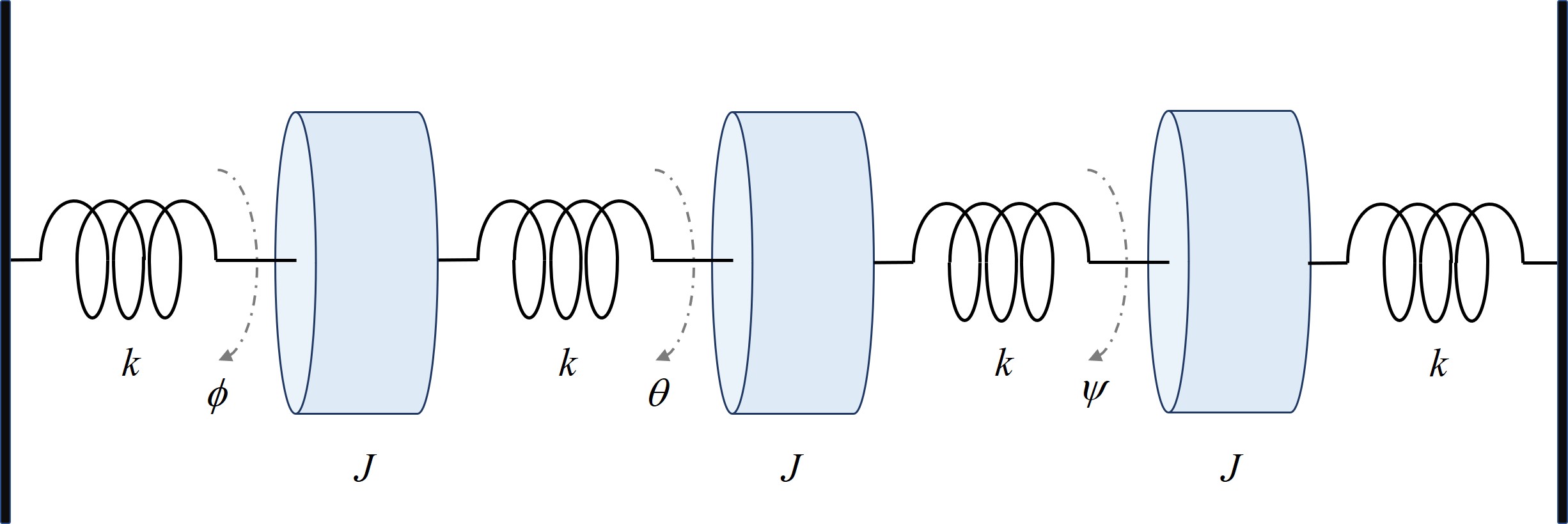}\\
  \caption{The three-inertia system.}\label{Fig-System}
\end{figure}

Consider a three-inertia system as shown in Figure \ref{Fig-System},
which is monitored by three separate sensors.
Denote each of the inertia's angle by $\phi, \theta$, and $\psi$, respectively,
and suppose that the sensors' measurements are
$y_{1}=\phi+\psi$, $y_{2}=\theta$, and $y_{3}=\psi-\phi$, respectively.
Then, with the state $x:=\text{col}\left(\phi, \dot{\phi}, \theta, \dot{\theta}, \psi, \dot{\psi}\right)$,
this system is in the form of \eqref{eq-x-Ax} and \eqref{eq-y-i} with
\begin{align*}
    A &= \left[
           \begin{array}{cccccc}
             0 & 1 & 0 & 0 & 0 & 0 \\
             -\frac{2k}{J} & 0 & \frac{k}{J} & 0 & 0 & 0 \\
             0 & 0 & 0 & 1 & 0 & 0 \\
             \frac{k}{J} & 0 & -\frac{2k}{J} & 0 & \frac{k}{J} & 0 \\
             0 & 0 & 0 & 0 & 0 & 1 \\
             0 & 0 & \frac{k}{J} & 0 & -\frac{2k}{J} & 0 \\
           \end{array}
         \right] \notag\\
    C&=\left[
      \begin{array}{c}
        C_{1} \\
        C_{2} \\
        C_{3} \\
      \end{array}
    \right]=\left[
              \begin{array}{cccccc}
                1 & 0 & 0 & 0 & 1 & 0 \\
                0 & 0 & 1 & 0 & 0 & 0 \\
                -1 & 0 & 0 & 0 & 1 & 0 \\
              \end{array}
            \right]
\end{align*}
where $k$ is the torsional stiffness and $J$ is the moment of the inertia.
It can be verified that for any positive $k$ and $J$,
the matrix $A$ is neutrally stable and the pair $(C,A)$ is observable.
Thus, Assumptions \ref{Ass-A-neutral} and \ref{Ass-observable} are satisfied.
Nevertheless, none of the pairs $(C_{i},A), i=1,2,3$, is observable.
In particular, $\text{rank}(\mathcal{O}_{1})=\text{rank}(\mathcal{O}_{2})=4$
and $\text{rank}(\mathcal{O}_{3})=2$.

\begin{figure}
  \centering
  \includegraphics[scale=0.3]{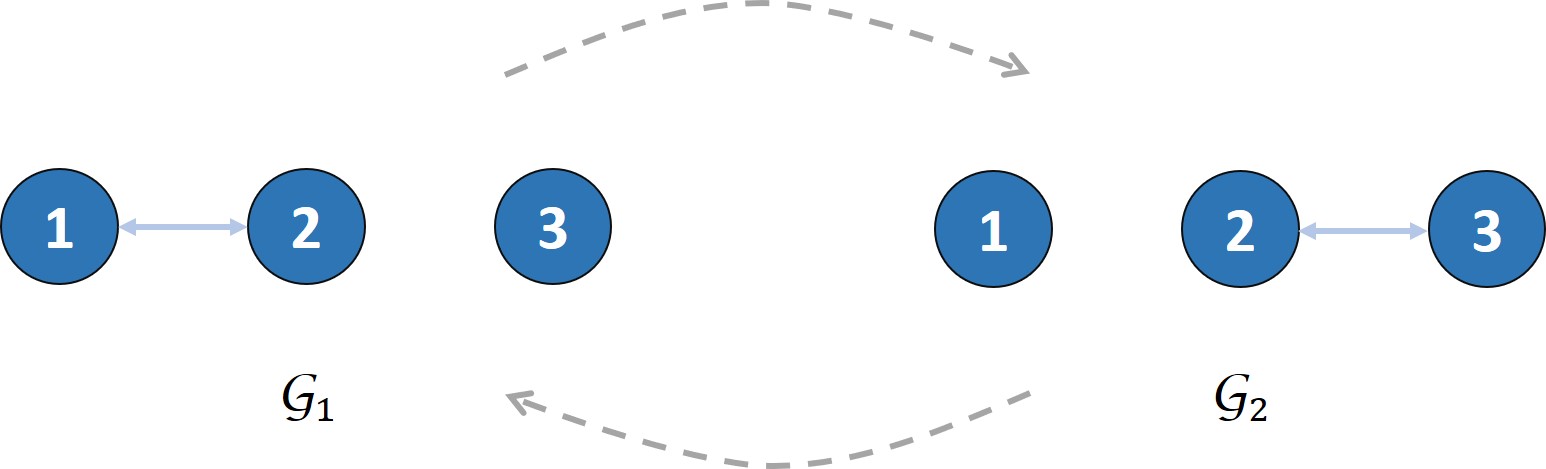}\\
  \caption{The switching graph $\mathcal{G}_{\sigma(t)}$.}\label{Fig-Graph}
\end{figure}

Suppose the three sensors transmit information over a communication network
described by the switching graph in Figure \ref{Fig-Graph},
which is dictated by the following switching signal:
\begin{equation*}
    \sigma(t)=
    \begin{cases}
   1, & \text{if}\quad  s T_{\text{c}} \le t < \left(s+\frac{1}{3}\right) T_{\text{c}}   \\
   2, & \text{if}\quad  \left(s+\frac{1}{3}\right)T_{\text{c}} \le t < (s+1) T_{\text{c}}
    \end{cases}
\end{equation*}
where $s=0,1,2,\ldots$. Then, clearly, Assumption \ref{Ass-jointly-connected} is verified.

\begin{figure}
\centering
  \includegraphics[trim=150 0 0 0, scale=0.12]{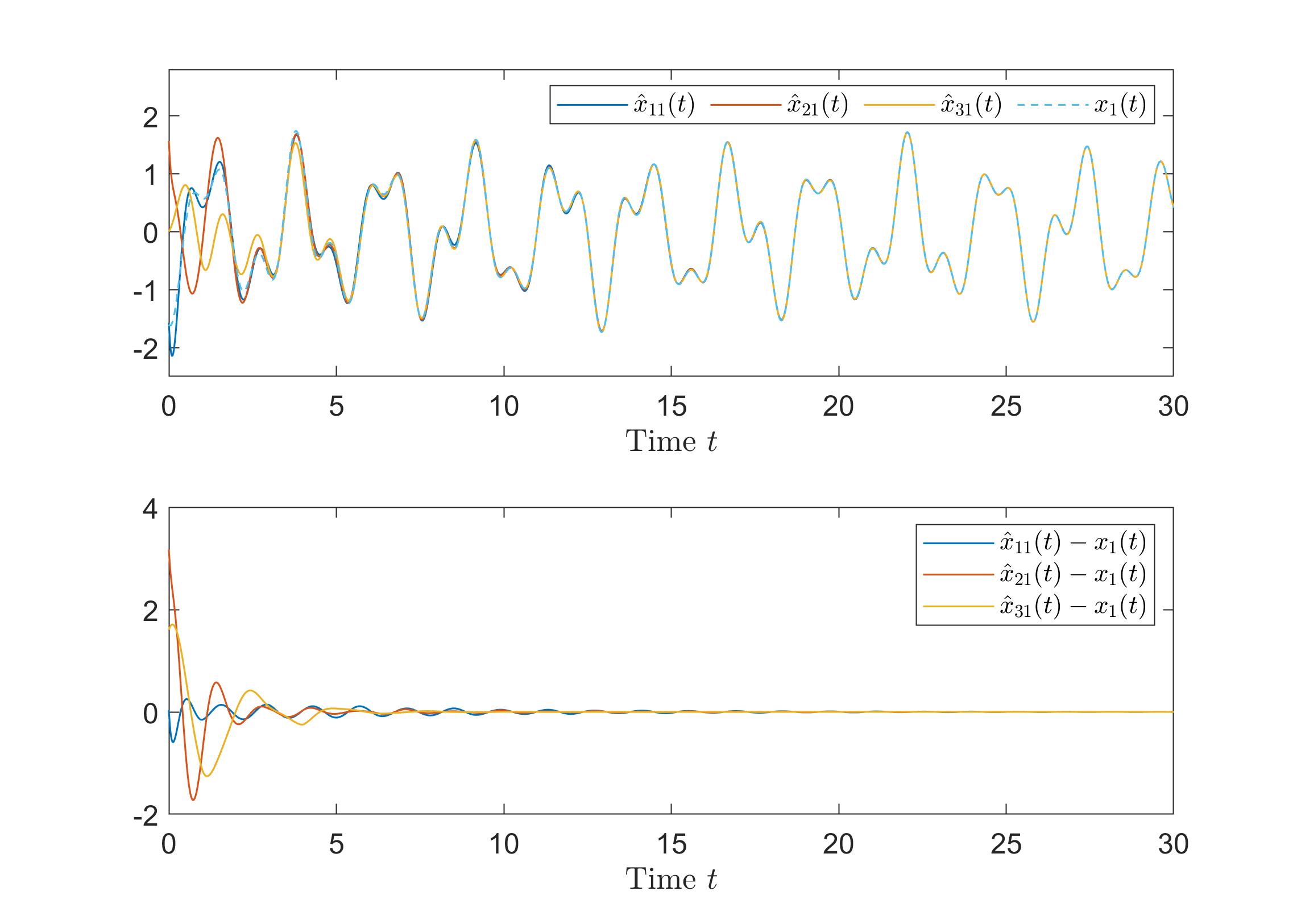}\\
  \caption{States $\hat{x}_{i1}(t), x_{1}(t)$ and estimation errors $\hat{x}_{i1}(t)-x_{1}(t),i=1,2,3$.}\label{Fig-Estimate_1}
\end{figure}

\begin{figure}
\centering
  \includegraphics[trim=150 0 0 0, scale=0.12]{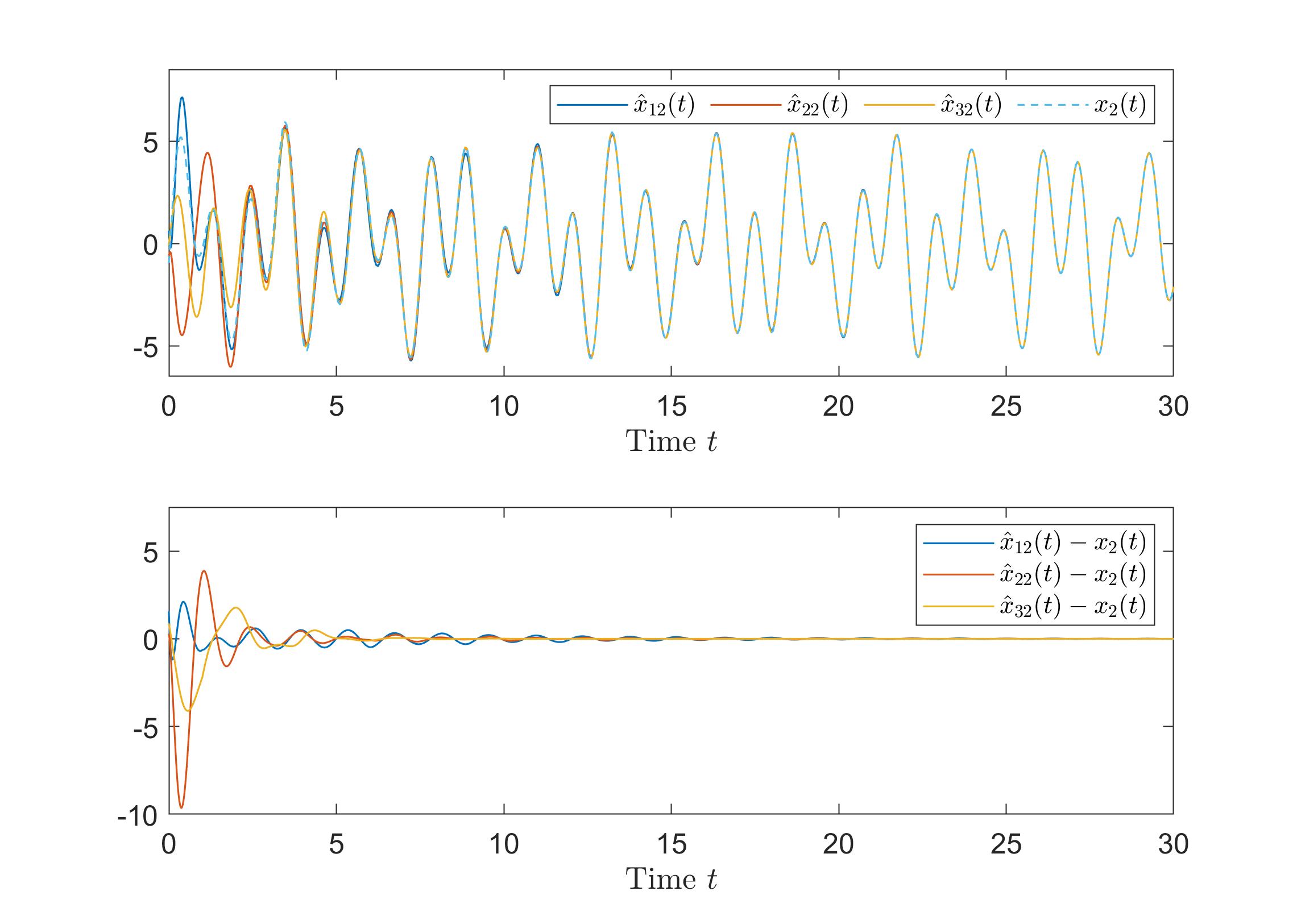}\\
  \caption{States $\hat{x}_{i2}(t), x_{2}(t)$ and estimation errors $\hat{x}_{i2}(t)-x_{2}(t),i=1,2,3$.}\label{Fig-Estimate_2}
\end{figure}

\begin{figure}
\centering
  \includegraphics[trim=150 0 0 0, scale=0.12]{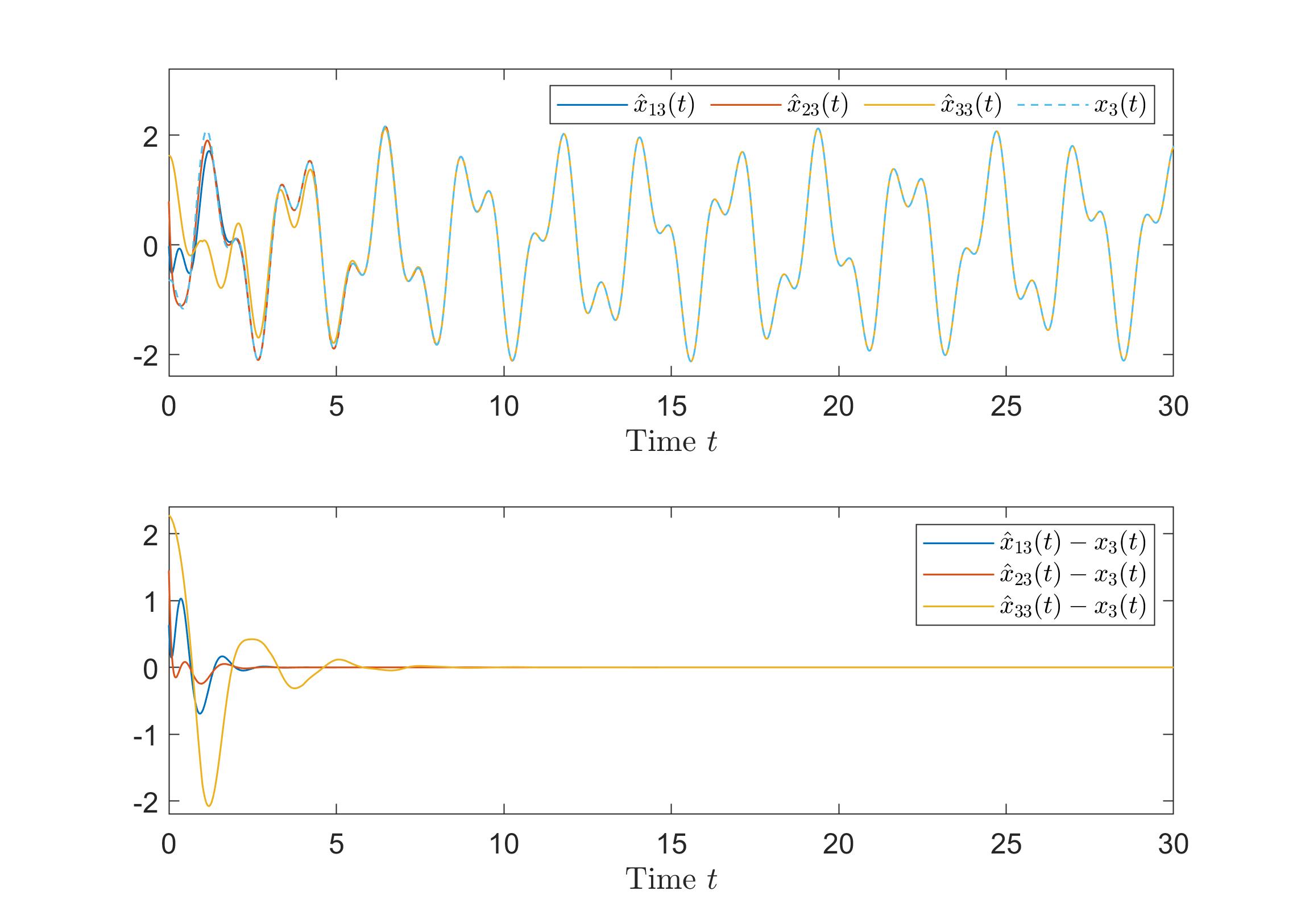}\\
  \caption{States $\hat{x}_{i3}(t), x_{3}(t)$ and estimation errors $\hat{x}_{i3}(t)-x_{3}(t),i=1,2,3$.}\label{Fig-Estimate_3}
\end{figure}

\begin{figure}
\centering
  \includegraphics[trim=150 0 0 0, scale=0.12]{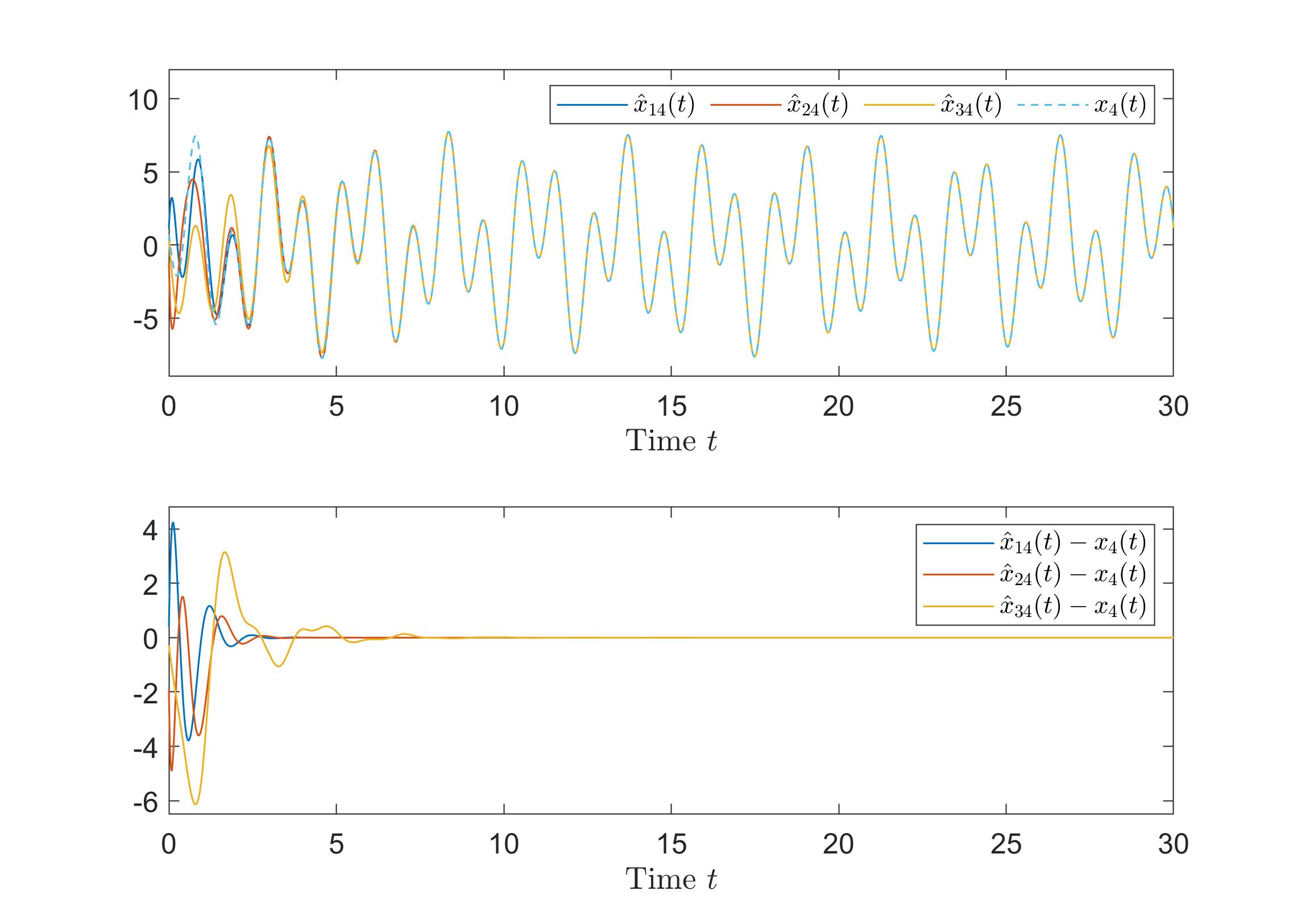}\\
  \caption{States $\hat{x}_{i4}(t), x_{4}(t)$ and estimation errors $\hat{x}_{i4}(t)-x_{4}(t),i=1,2,3$.}\label{Fig-Estimate_4}
\end{figure}

\begin{figure}
\centering
  \includegraphics[trim=150 0 0 0, scale=0.12]{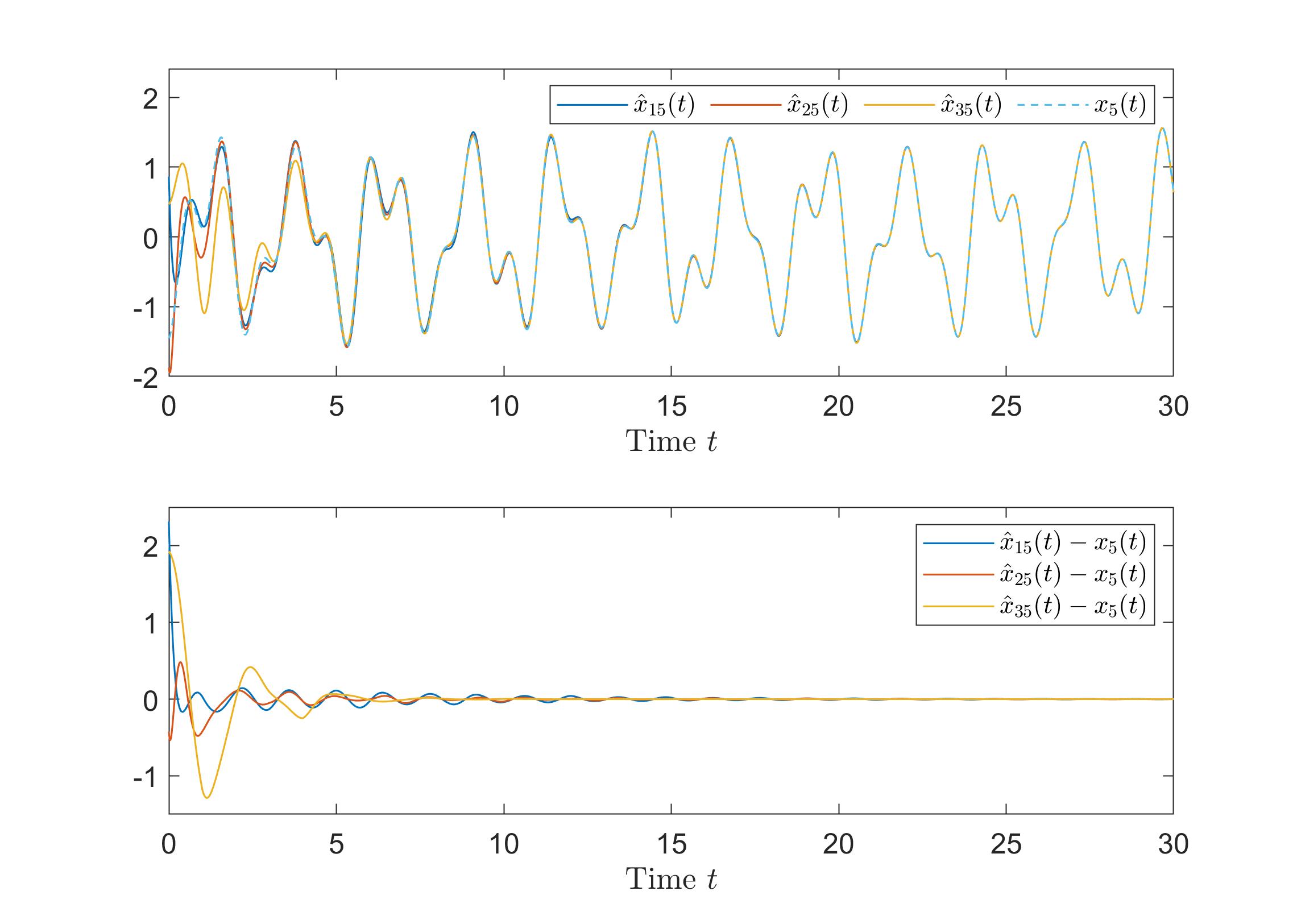}\\
  \caption{States $\hat{x}_{i5}(t), x_{5}(t)$ and estimation errors $\hat{x}_{i5}(t)-x_{5}(t),i=1,2,3$.}\label{Fig-Estimate_5}
\end{figure}

\begin{figure}
\centering
  \includegraphics[trim=150 0 0 0, scale=0.12]{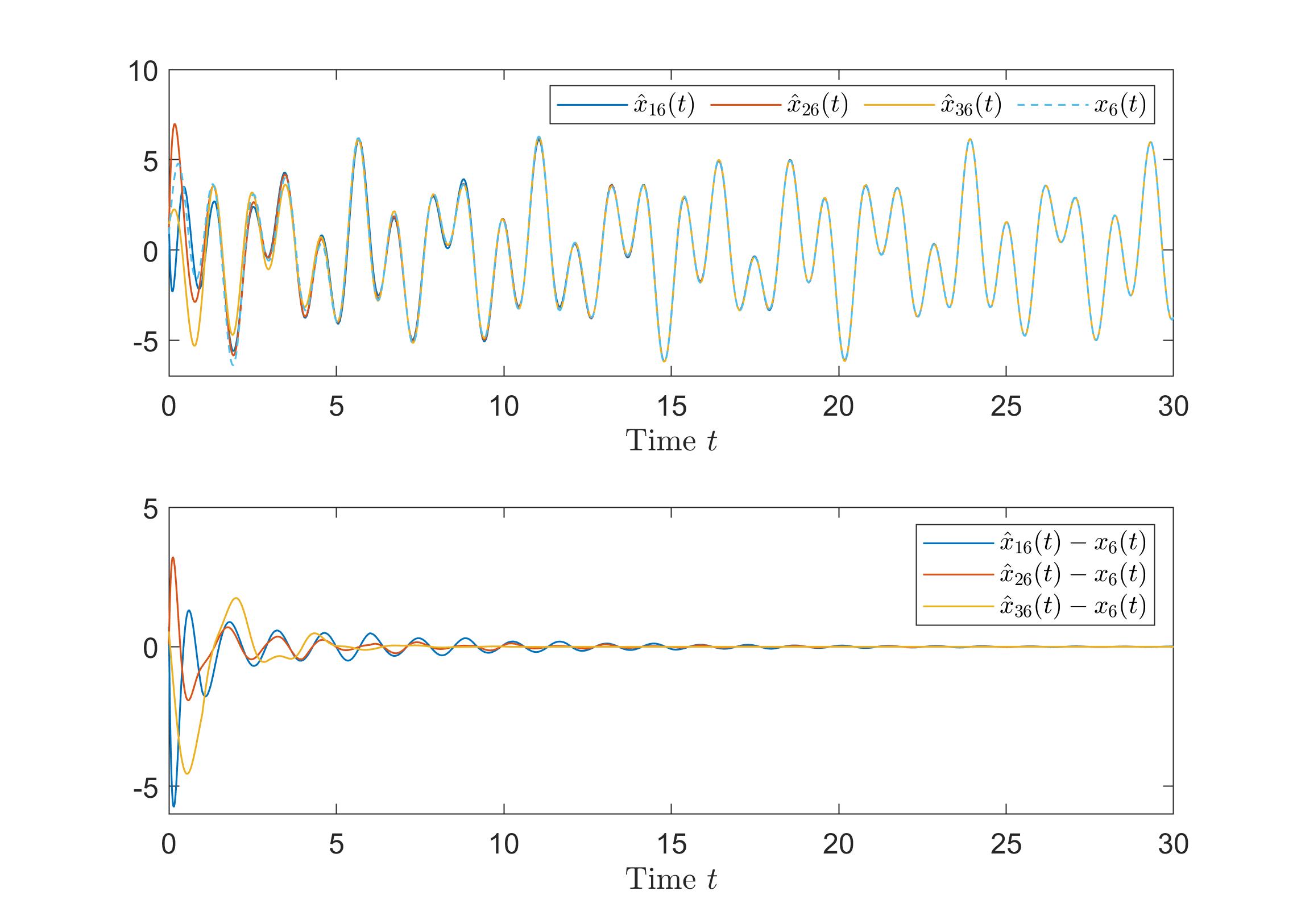}\\
  \caption{States $\hat{x}_{i6}(t), x_{6}(t)$ and  estimation errors $\hat{x}_{i6}(t)-x_{6}(t),i=1,2,3$.}\label{Fig-Estimate_6}
\end{figure}

Thus, by Theorem \ref{Theorem}, we can solve Problem \ref{Problem}
by designing, for $i=1,2,3$, a local observer of the form \eqref{eq-hat-x-i-specific}.
Suppose $\frac{k}{J}=10\, \mathrm{N} \cdot (\mathrm{m \cdot kg \cdot rad})^{-1}$.
Following the design procedure sketched in Remark \ref{Remark-design-procedure},
we choose $L_{i\text{o}}, i=1,2$, such that the eigenvalues of $(A_{i\text{o}}-L_{i\text{o}}C_{i\text{o}})$ are
placed at $\{-2\pm j 5, -5\pm j 2\}$,
and $L_{3\text{o}}$ to place the eigenvalues of $(A_{3\text{o}}-L_{3\text{o}}C_{3\text{o}})$ at $\{-2,-5\}$.
The performance of the local observers is simulated with
$T_{\text{c}}=3$, $\gamma_{i}=1$, $a_{ij}(t)=1$ if $(j,i) \in \mathcal{E}_{\sigma(t)}$, $i,j=1,2,3$, and
randomly generated initial conditions.
Figures \ref{Fig-Estimate_1} to \ref{Fig-Estimate_6} show, respectively,
each component of the state of the local observers and the system,
together with each component of the estimation errors.

\section{Conclusion}\label{Section-Conclusion}
In this paper, we have presented an exponential convergence result on the distributed state estimation problem for linear systems over jointly connected switching networks.
The main result
offers two advantages over the existing one.
First, the exponential convergence leads to the guaranteed convergence rate, which is much  desired in practice.
Second, since the error system is uniformly asymptotically stable,
it is also totally stable and hence is able to withstand small disturbances.
These two advantages are achieved by establishing exponential stability for two classes of linear switched systems,
which may have some other applications.
A restriction of the current result is that it only applies to marginally stable linear systems.
Thus, it would be interesting to further consider removing or relaxing this restriction,
so as to accommodate more general linear systems.
It would also be interesting to consider extending the result to directed switching networks.

\section*{Appendix A: Notation on graph}
A graph $\mathcal{G}:=\left(\mathcal{V}, \mathcal{E}\right)$ consists of a finite node set $\mathcal{V}:=\{1, 2, \ldots, N\}$
and an edge set $\mathcal{E} \subseteq \mathcal{V} \times \mathcal{V}$.
An edge of $\mathcal{E}$ from node $j$ to node $i$, $j \ne i$,
is denoted by $(j,i)$, and node $j$ is called a neighbor of node $i$.
Then, $\mathcal{N}_{i}:=\left\{j \in \mathcal{V} : (j,i) \in \mathcal{E} \right\}$ is called the neighbor set of node $i$.
The edge $(i,j)$ is called undirected if $(i,j) \in \mathcal{E}$ implies $(j,i) \in \mathcal{E}$.
The graph $\mathcal{G}$ is called undirected if every edge in $\mathcal{E}$ is undirected.
If the graph contains a set of edges of the form $\{(i_{1},i_{2})$, $(i_{2},i_{3})$,
$\ldots$, $(i_{k-1},i_{k})\}$, then this set is called a path from node $i_{1}$ to node $i_{k}$,
and node $i_{k}$ is said to be reachable from node $i_{1}$.
A graph is called strongly connected if there exists a path between any two distinct nodes.
An undirected and strongly connected graph is called connected.

The weighted adjacency matrix of a graph $\mathcal{G}$ is a nonnegative matrix
$\mathcal{A}:=[a_{ij}]_{i,j=1}^{N} \in \RR^{N \times N}$, where
$a_{ii}=0$ and,
for $i\ne j, a_{ij}>0$ if and only if $ (j,i)\in \mathcal{E}$.
Then, the Laplacian $\mathcal{L}:=[l_{ij}]_{i,j=1}^{N} \in \RR^{N \times N}$ of the graph $\mathcal{G}$
can be defined from $\mathcal{A}$ with $l_{ii}=\sum_{j=1}^{N}a_{ij}$ and,
for $i \ne j, l_{ij}=-a_{ij}$.
Moreover, $\mathcal{L}$ is symmetric and positive semi-definite if and only if the graph $\mathcal{G}$ is undirected \cite{Godsil01}.

Given the switching signal $\sigma: [0,\infty) \mapsto \mathcal{P}=\{1,2,\ldots,n_{0}\}$ and $n_{0}$ graphs
$\mathcal{G}_{p}=(\mathcal{V}, \mathcal{E}_{p})$, $p=1,2,\ldots,n_{0}$,
each with the corresponding weighted adjacency matrix denoted by $\mathcal{A}_{p}$
and the Laplacian by $\mathcal{L}_{p}$, $p=1,2, \ldots, n_{0}$,
we call the time-varying graph $\mathcal{G}_{\sigma(t)}=\left(\mathcal{V},\mathcal{E}_{\sigma(t)}\right)$ a switching graph and denote
its weighted adjacency matrix by $\mathcal{A}_{\sigma(t)}$,
and its Laplacian by $\mathcal{L}_{\sigma(t)}$.
Finally, the graph $\mathcal{G}=\left(\mathcal{V}, \mathcal{E}\right)$ with
$\mathcal{E}= \bigcup_{p=1}^{r} \mathcal{E}_{p}$ is called the union of the graphs $\mathcal{G}_{p}, p=1,2,\ldots,r$,
and is denoted by $\mathcal{G}=\bigcup_{p=1}^{r} \mathcal{G}_{p}$.

\section*{Appendix B: Proof of Lemma \ref{Lemma-eta-stable}}
To prove Lemma \ref{Lemma-eta-stable}, we first need to show the following result.

\begin{Lemma}\label{Lemma-UCO}
Consider the following system:
\begin{equation}\label{eq-bar-eta-system}
    \dot{\bar{\eta}}(t)=0, \qquad \bar{y}(t)=\left(\mathcal{L}_{\sigma(t)}^{\frac{1}{2}} \otimes I_{n}\right)U \bar{\eta}(t), \qquad t\ge 0
\end{equation}
where $\bar{\eta}(t)\in \RR^{\nu}$ is the state and $\bar{y}(t) \in \RR^{Nn}$ is the output.
Under Assumptions \ref{Ass-observable} and \ref{Ass-jointly-connected},
  system \eqref{eq-bar-eta-system} is UCO, i.e.,
  there exist $T_{\text{o}}>0$ and $0<\bar{\alpha}_{1}\le \bar{\alpha}_{2}$,
  such that the observability Gramian
  $G_{\bar{\eta}}(t^{*},t^{*}+T_{\text{o}})$ of system \eqref{eq-bar-eta-system} satisfies
$
    \bar{\alpha}_{1}I_{\nu} \le G_{\bar{\eta}}(t^{*},t^{*}+T_{\text{o}}) \le \bar{\alpha}_{2}I_{\nu}, \forall\ t^{*}\ge 0.
$
\end{Lemma}

\begin{Proof}
  Let $T_{\text{o}}\ge 2T_{\text{c}}$.
Then, under Assumption \ref{Ass-jointly-connected},  for any fixed $t^{*}\ge 0$,
there exists some switching instant $t_{j_{k}}$ that satisfies
$
    t^{*} \le t_{j_{k}} < t_{j_{k+1}} < t^{*} +2 T_{\text{c}}.
$
Specifically, the switching instants within $[t_{j_{k}}, t_{j_{k+1}})$ can be listed as
$\{t_{j_{k}}, t_{j_{k}+1}, t_{j_{k}+2}, \ldots, t_{j_{k+1}-1}\}$.
Thus, we have
\begin{align}\label{eq-Gramian-bar-eta-ge}
  &\quad\  G_{\bar{\eta}}(t^{*},t^{*}+T_{\text{o}})  \notag \\
  &=\int_{t^{*}}^{t^{*}+T_{\text{o}}} U^{T}\left(\mathcal{L}^{\frac{1}{2}}_{\sigma(t)}\otimes I_{n} \right)^{T}\left(\mathcal{L}^{\frac{1}{2}}_{\sigma(t)}\otimes I_{n}\right)U  d\, t \notag \\
  & \ge\int_{t^{*}}^{t^{*}+2 T_{\text{c}}}U^{T}\left(\mathcal{L}_{\sigma(t)}\otimes I_{n} \right)U  d\, t  \notag \\
  &\ge \int_{t_{j_{k}}}^{t_{j_{k+1}}} U^{T}\left(\mathcal{L}_{\sigma(t)}\otimes I_{n} \right)U d\, t  \notag \\
  &= \left(\int_{t_{j_{k}}}^{t_{j_{k}+1}}+\cdots+\int_{t_{j_{k+1}-1}}^{t_{j_{k+1}}}\right)
U^{T}\left(\mathcal{L}_{\sigma(t)}\otimes I_{n} \right)U d\, t \notag\\
& \ge \tau \cdot U^{T}\left(\sum_{r=j_{k}}^{j_{k+1}-1}\mathcal{L}_{\sigma(t_{r})}\otimes I_{n} \right)U .
\end{align}

Next, we show that under Assumptions \ref{Ass-observable} and \ref{Ass-jointly-connected},
there exists some $\hat{\alpha}_{1}>0$ such that
  \begin{equation}\label{eq-U-L-U-pd}
    U^{T}\left(\sum_{r=j_{k}}^{j_{k+1}-1}\mathcal{L}_{\sigma(t_{r})}\otimes I_{n} \right)U \ge \hat{\alpha}_{1} I_{\nu},
    \quad    \forall \, k=0,1,2,\ldots.
  \end{equation}
If so, combining \eqref{eq-Gramian-bar-eta-ge} and \eqref{eq-U-L-U-pd} yields
$ G_{\bar{\eta}}(t^{*},t^{*}+T_{\text{o}}) \ge \bar{\alpha}_{1}I_{\nu}, \forall\ t^{*}\ge 0$
with $\bar{\alpha}_{1}:=\tau \hat{\alpha}_{1}$.

Clearly, $U^{T}\left(\sum_{r=j_{k}}^{j_{k+1}-1}\mathcal{L}_{\sigma(t_{r})}\otimes I_{n} \right)U$
  is symmetric and positive semi-definite for any $k=0,1,2,\ldots$.
  Suppose for some
  $v :=\text{col}(v_{1},\ldots, v_{N})$ with $v_{i}\in \RR^{\nu_{i}}, i=1,\ldots,N$,
  the following holds:
  \begin{equation}\label{eq-L-U-x-0}
    \left(\sum_{r=j_{k}}^{j_{k+1}-1}\mathcal{L}_{\sigma(t_{r})}\otimes I_{n} \right)U v=0.
  \end{equation}
 Under Assumption \ref{Ass-jointly-connected}, by Remark \ref{Remark-jointly-connected},
  the null space of the matrix $\left(\sum_{r=j_{k}}^{j_{k+1}-1}\mathcal{L}_{\sigma(t_{r})}\otimes I_{n} \right)$
  is spanned by the columns of the matrix $\mathbf{1}_{N}\otimes I_{n}$.
 It then follows from \eqref{eq-U-block-def} and \eqref{eq-L-U-x-0} that
  \begin{equation}\label{eq-y-U-i-x-i}
    z:=U_{1}v_{1}=U_{2}v_{2}=\cdots=U_{N}v_{N}
  \end{equation}
 for some $z \in \RR^{n}$.
 By the definition of $U_{i}, i=1,\ldots,N$, in \eqref{eq-U-i-def},
 we see that $z \in \text{ker}(\mathcal{O}_{i}), i=1,\ldots,N$.
 Thus, $z\in \bigcap_{i=1}^{N}\text{ker}(\mathcal{O}_{i})$, which,
 under Assumption \ref{Ass-observable}, implies that $z=0$.
 Since $U_{i}, i=1,\ldots,N$, are of full column rank,
 from \eqref{eq-y-U-i-x-i}, we have $v_{1}=0, \ldots, v_{N}=0$, and hence $v=0$.
 This shows that for any $k=0,1,2,\ldots$,
 $U^{T}\left(\sum_{r=j_{k}}^{j_{k+1}-1}\mathcal{L}_{\sigma(t_{r})}\otimes I_{n} \right)U$ is positive definite.
 By further noting that $\sigma(t)$ only takes on finitely many values
 and that $t_{j_{k+1}}-t_{j_{k}}, k=0,1,2,\ldots$, are uniformly bounded by the finite $T_{\text{c}}$,
 we conclude that there exists some $\bar{\alpha}_{1}>0$ such that \eqref{eq-U-L-U-pd} holds.

Finally, the existence of $\bar{\alpha}_{2}>0$ satisfying
$G_{\bar{\eta}}(t^{*},t^{*}+T_{\text{o}}) \le \bar{\alpha}_{2}I_{\nu}, \forall\ t^{*}\ge 0$
is obvious, since both $T_{\text{o}}$ and the range of $\sigma(t)$ are finite.
The overall proof is thus complete.
\end{Proof}

\begin{Remark}
The  positive definiteness of the matrix $U^{T}\left(\sum_{r=j_{k}}^{j_{k+1}-1}\mathcal{L}_{\sigma(t_{r})}\otimes I_{n} \right)U$ under Assumptions \ref{Ass-observable} and \ref{Ass-jointly-connected} was also asserted in Lemma 3.3 of \cite{ZhangLu21ACC} without the proof.
\end{Remark}

\medskip \noindent
\emph{Proof of Lemma \ref{Lemma-eta-stable}:}
First of all,
let us specify an output $y(t)\in \RR^{Nn}$ for system \eqref{eq-eta-system} as follows:
\begin{equation}\label{eq-y-specified-output}
    y(t)=\left(\mathcal{L}_{\sigma(t)}^{\frac{1}{2}} \otimes I_{n}\right)U\eta(t).
\end{equation}
Then, we observe that system \eqref{eq-eta-system} with the output \eqref{eq-y-specified-output} is in the following form:
\begin{equation*}
\dot{\eta}(t) = \left(\mathbf{A}-\mathbf{F}(t)\mathbf{C}(t)\right)\eta(t),
\qquad y(t)=\mathbf{C}(t)\eta(t)
\end{equation*}
with $\mathbf{A}=\mathbf{0}$ and
\begin{equation*}
   \mathbf{F}(t):=Q U^{T}\left(\mathcal{L}^{\frac{1}{2}}_{\sigma(t)} \otimes I_{n}\right), \quad \mathbf{C}(t):=\left(\mathcal{L}_{\sigma(t)}^{\frac{1}{2}} \otimes I_{n}\right)U
\end{equation*}
which are uniformly bounded over $[0,\infty)$.
Thus, by Lemma 1 of \cite{Anderson77},
system \eqref{eq-eta-system} with the output \eqref{eq-y-specified-output} is UCO if and only if system \eqref{eq-bar-eta-system} is UCO.
As a result of Lemma \ref{Lemma-UCO}, under Assumptions \ref{Ass-observable} and \ref{Ass-jointly-connected},
there exist  $0<\alpha_{1}\le \alpha_{2}$,
such that the observability Gramian
  $G_{\eta}(t^{*},t^{*}+T_{\text{o}})$ of system \eqref{eq-eta-system} with the output \eqref{eq-y-specified-output} satisfies
\begin{equation}\label{eq-Gramian-inequality}
    \alpha_{1}I_{\nu} \le G_{\eta}(t^{*},t^{*}+T_{\text{o}}) \le \alpha_{2}I_{\nu},\quad  \forall\ t^{*}\ge 0.
\end{equation}

Define $V(\eta)=\frac{1}{2}\eta^{T}Q^{-1}\eta$.
Then, the time derivative of $V(\eta(t))$ along the trajectory of system \eqref{eq-eta-system} satisfies
\begin{align}\label{eq-V-dot}
    \dot{V}(\eta(t))&= - \eta(t)^{T} U^{T}
    \left(\mathcal{L}_{\sigma(t)} \otimes I_{n}\right)U\eta(t)  \notag\\
    &=- \eta(t)^{T}\mathbf{C}(t)^{T}\mathbf{C}(t)\eta(t) \le 0.
\end{align}
Hence, we have $V(\eta(t))\le V(\eta(t^{*})),\forall\, t\ge t^{*}$ and
\begin{equation*}
    \|\eta(t)\| \le  \sqrt{\frac{\lambda_{\max}(Q)}{\lambda_{\min}(Q)}} \|\eta(t^{*})\|, \qquad \forall\, t\ge t^{*}.
\end{equation*}
Thus, system \eqref{eq-eta-system} is uniformly stable.

Next, denote the state transition matrix of system \eqref{eq-eta-system}
by $\Phi_{\eta}(t,t^{*}), t\ge t^{*} \ge 0$.
Then, it follows from \eqref{eq-V-dot} that
\begin{align*}
   &\quad\ V(\eta(t^{*}+T_{\text{o}}))-V(\eta(t^{*})) \notag\\
   &=- \int_{t^{*}}^{t^{*}+T_{\text{o}}} \eta(t)^{T}\mathbf{C}(t)^{T}\mathbf{C}(t)\eta(t)\, d \, t\notag\\
   &=- \eta(t^{*})^{T}\Bigg( \int_{t^{*}}^{t^{*}+T_{\text{o}}} \Phi_{\eta}(t,t^{*})^{T}\mathbf{C}(t)^{T} \notag\\
   &\qquad\qquad\qquad\qquad\qquad\qquad\quad \times   \mathbf{C}(t)\Phi_{\eta}(t,t^{*})\, d \, t \Bigg)\eta(t^{*}) \notag\\
   &=- \eta(t^{*})^{T} G_{\eta}(t^{*},t^{*}+T_{\text{o}})\eta(t^{*}).
\end{align*}
By further noting \eqref{eq-Gramian-inequality}, we obtain
\begin{align}\label{eq-V-t+T-o-t}
&\quad\ V(\eta(t^{*}+T_{\text{o}}))-V(\eta(t^{*})) \notag\\
& \le -\alpha_{1} \eta(t^{*})^{T} \eta(t^{*}) \notag \\
& \le -2\alpha_{1}\lambda_{\min}(Q) V(\eta(t^{*})) \le  -\rho V(\eta(t^{*}))
\end{align}
where $\rho$ is chosen to satisfy $0< \rho < \min\{2\alpha_{1} \lambda_{\min}(Q), 1\}$.
By rearranging the terms in \eqref{eq-V-t+T-o-t}, we have
\begin{equation}\label{eq-V-strict-decrease}
    V(\eta(t^{*}+T_{\text{o}}))\le (1-\rho)V(\eta(t^{*})).
\end{equation}

Now, given any $\delta>0$, there exists a positive integer $\ell$ such that
$
    (1-\rho)^{\ell} \le \delta^{2} \frac{\lambda_{\min}(Q)}{\lambda_{\max}(Q)}.
$
Then, by \eqref{eq-V-dot} and \eqref{eq-V-strict-decrease},
for any $t \ge t^{*}+\ell T_{\text{o}}$, we have
\begin{align*}
    V(\eta(t)) & \le V(\eta(t^{*}+\ell T_{\text{o}}))\\
    & \le (1-\rho)^{\ell}V(\eta(t^{*})) \le  \delta^{2} \frac{\lambda_{\min}(Q)}{\lambda_{\max}(Q)} V(\eta(t^{*}))
\end{align*}
which implies that
\begin{equation*}\label{}
    \|\eta(t)\| \le \delta \|\eta(t^{*})\|, \quad \forall\, t \ge t^{*}+\ell T_{\text{o}}.
\end{equation*}
Thus, system \eqref{eq-eta-system} is uniformly asymptotically stable.
Since, by Theorem 6.13 of \cite{Rugh96},
uniform asymptotical stability and exponential stability are equivalent for linear time-varying systems,
we conclude that system \eqref{eq-eta-system} is exponentially stable.
The proof is complete.   \hfill $\Box$

%%%%%%%%%%%%%%%%%%%%%%%%%%%%%%%%%%%%%%%%%%%%%%%%%%%%%%%%%%%%%%%%%%%%%%%%%%%%%%%%%%%%%%%%%%%%%%%%%%%%%%%%%%%%%%%%%%%%%%%%%%%%%%%%%%%%%%%%%%%%%%%%%%%%%%%%%%%%%%%%%%%%%%%%%%%%%%%%%%%%%%%%%%%%%%%%%%%%%%%%%%%%%%%

\end{document}